\documentclass[a4paper]{amsart}

\usepackage{amsmath}
\usepackage{amsfonts}
\usepackage{palatino}
\usepackage{a4wide}
\usepackage{amssymb}
\usepackage{amsthm}

\theoremstyle{plain}

\newtheorem{teo}{Theorem}[section]

\newtheorem{prop}[teo]{Proposition}
\newtheorem{cor}[teo]{Corollary}
\newtheorem{ackn}{Acknowledgments\!}

\theoremstyle{definition}

\theoremstyle{remark}

\newtheorem{rem}[teo]{Remark}

\numberwithin{equation}{section}

\def\SS{{{\mathbb S}}}
\def\NN{{{\mathbb N}}}

\def\RR{{\mathbb R}}

\def\HH{{\mathbb H}}
\def\RRR{{\mathrm R}}
\def\WWW{{\mathrm W}}
\def\AAA{{\mathrm A}}
\def\BBB{{\mathrm B}}
\def\DDD{{\mathrm D}}

\def\Ric{{\mathrm {Ric}}}

\def\SSS{{\mathrm S}}
\def\CCC{{\mathrm C}}
\def\Rm{{\mathrm {Rm}}}

\def\dt{\frac{\partial\,}{\partial t}}

\def\div{\operatornamewithlimits{div}\nolimits}

\title[The Evolution of the Weyl Tensor under the Ricci Flow]{The Evolution of the Weyl Tensor under the Ricci Flow}
\date{\today}

\author[Giovanni Catino]{Giovanni Catino}
\address[Giovanni Catino]{SISSA -- International School for Advanced
  Studies, Via Bonomea 265, Trieste, Italy, 34136}
\email[G. Catino]{catino@sissa.it}

\author[Carlo Mantegazza]{Carlo Mantegazza}
\address[Carlo Mantegazza]{Scuola Normale Superiore di Pisa, P.za Cavalieri 7, Pisa, Italy, 56126}
\email[C. Mantegazza]{c.mantegazza@sns.it}

\date{\today}

\begin{document}

\begin{abstract} We compute the evolution equation of the Weyl tensor
  under the Ricci flow of a Riemannian manifold and we discuss some consequences
  for the classification of locally conformally flat Ricci solitons.
\end{abstract}

\maketitle
\tableofcontents

\section{The Evolution Equation of the Weyl Tensor}

The Riemann curvature
operator of a Riemannian manifold $(M^n,g)$ is defined 
as in~\cite{gahula} by
$$
\mathrm{Riem}(X,Y)Z=\nabla_{Y}\nabla_{X}Z-\nabla_{X}\nabla_{Y}Z+\nabla_{[X,Y]}Z\,.
$$ 
In a local coordinate system the components of the $(3,1)$--Riemann 
curvature tensor are given by
$\RRR^{l}_{ijk}\tfrac{\partial}{\partial
  x^{l}}=\mathrm{Riem}\big(\tfrac{\partial}{\partial
  x^{i}},\tfrac{\partial}{\partial
  x^{j}}\big)\tfrac{\partial}{\partial x^{k}}$ and we denote by
$\RRR_{ijkl}=g_{lm}\RRR^{m}_{ijk}$ its $(4,0)$--version.\\

\medskip

{\em In all the paper the Einstein convention of summing over the repeated 
indices will be adopted.}

\medskip

With this choice, for the sphere $\SS^n$ we have
${\mathrm{Riem}}(v,w,v,w)=\RRR_{ijkl}v^iw^jv^kw^l>0$.

The Ricci tensor is obtained by the contraction 
$\RRR_{ik}=g^{jl}\RRR_{ijkl}$ and $\RRR=g^{ik}\RRR_{ik}$ will 
denote the scalar curvature.

The so called Weyl tensor is then 
defined by the following decomposition formula (see~\cite[Chapter~3,
Section~K]{gahula}) in dimension $n\geq 3$,
\begin{align*}
\WWW_{ijkl}=&\,\RRR_{ijkl}+\frac{\RRR}{(n-1)(n-2)}(g_{ik}g_{jl}-g_{il}g_{jk})
- \frac{1}{n-2}(\RRR_{ik}g_{jl}-\RRR_{il}g_{jk}
+\RRR_{jl}g_{ik}-\RRR_{jk}g_{il})\\
=&\,\RRR_{ijkl}+\AAA_{ijkl}+\BBB_{ijkl}\,,
\end{align*}
where we introduced the tensors
$$
\AAA_{ijkl}=\frac{\RRR}{(n-1)(n-2)}(g_{ik}g_{jl}-g_{il}g_{jk})
$$
and
$$
\BBB_{ijkl}=-\frac{1}{n-2}(\RRR_{ik}g_{jl}-\RRR_{il}g_{jk}
+\RRR_{jl}g_{ik}-\RRR_{jk}g_{il})\,.
$$
The Weyl tensor satisfies all the symmetries of the curvature tensor
and all its traces with the metric are zero, 
as it can be easily seen by the above formula.\\
In dimension three $\WWW$ is identically zero for every Riemannian
manifold $(M^3,g)$, 
it becomes relevant instead when $n\geq 4$ since its
nullity is a condition equivalent for $(M^n,g)$ to be {\em locally
  conformally flat}, that is, around every point $p\in M^n$ there is a conformal
deformation $\widetilde{g}_{ij}=e^fg_{ij}$ of the original metric $g$,
such that the new metric is flat, 
namely, the Riemann tensor associated to
$\widetilde{g}$ is zero in $U_p$ (here $f:U_p\to\RR$ is a smooth function defined in a open
neighborhood $U_p$ of $p$).

We suppose now that $(M^n,g(t))$ is a Ricci flow in some time interval,
that is, the time--dependent metric $g(t)$ satisfies 
$$
\dt g_{ij} =-2 \RRR_{ij}\,.
$$
We have then the following evolution equations for the curvature 
(see for instance~\cite{hamilton1}),
\begin{align}
&\dt\RRR=\Delta\RRR+2\vert\Ric\vert^2\,\nonumber\\
&\dt\RRR_{ij}=\Delta\RRR_{ij}+2\RRR^{kl}\RRR_{kilj}
-2g^{pq}\RRR_{ip}\RRR_{jq}\,,\nonumber\\
&\dt\RRR_{ijkl}=\Delta\RRR_{ijkl}+2(\CCC_{ijkl}-\CCC_{ijlk}+\CCC_{ikjl}-\CCC_{iljk})\label{riemevol}\\
&\phantom{\dt\RRR_{ijkl}=}-g^{pq}(\RRR_{ip}\RRR_{qjkl}+\RRR_{jp}\RRR_{iqkl}+\RRR_{kp}\RRR_{ijql}+\RRR_{lp}\RRR_{ijkq})\nonumber\,,
\end{align}
where $\CCC_{ijkl}=g^{pq}g^{rs}\RRR_{pijr}\RRR_{slkq}$.

\medskip

{\em All the computations which follow 
will be done in a fixed local frame, not in
  a moving frame.}

\medskip

The goal of this section is to work out the evolution equation under
the Ricci flow of the Weyl tensor $\WWW_{ijkl}$. In the next sections
we will see the geometric consequences of the 
assumption that a manifold evolving by the Ricci flow is locally
conformally flat at every time. In particular, 
we will be able to classify the so called Ricci solitons under the
hypothesis of locally conformally flatness.

\bigskip

Since $\WWW_{ijkl}=\RRR_{ijkl}+\AAA_{ijkl}+\BBB_{ijkl}$ and we already
have the evolution equation~\eqref{riemevol} for $\RRR_{ijkl}$, we 
start differentiating in time the tensors $\AAA_{ijkl}$ and $\BBB_{ijkl}$
\begin{align*}
\dt\AAA_{ijkl}=&\,\frac{\Delta\RRR+2\vert\Ric\vert^2}{(n-1)(n-2)}(g_{ik}g_{jl}-g_{il}g_{jk})\\
&\,+\frac{\RRR}{(n-1)(n-2)}
(-2\RRR_{ik}g_{jl}-2\RRR_{jl}g_{ik}+2\RRR_{il}g_{jk}+2\RRR_{jk}g_{il})\\
=&\,\Delta\AAA_{ijkl}+\frac{2\vert\Ric\vert^2}{(n-1)(n-2)}(g_{ik}g_{jl}-g_{il}g_{jk})
+\frac{2\RRR}{n-1}\BBB_{ijkl}
\end{align*}
and
\begin{align*}
\dt\BBB_{ijkl}=&\,-\frac{1}{n-2}\Bigl((\Delta\RRR_{ik}+2\RRR^{pq}\RRR_{piqk}
-2g^{pq}\RRR_{ip}\RRR_{kq})g_{jl}\\
&\,-(\Delta\RRR_{il}+2\RRR^{pq}\RRR_{piql}-2g^{pq}\RRR_{ip}\RRR_{lq})g_{jk}\\
&\,+(\Delta\RRR_{jl}+2\RRR^{pq}\RRR_{pjql}-2g^{pq}\RRR_{jp}\RRR_{lq})g_{ik}\\
&\,-(\Delta\RRR_{jk}+2\RRR^{pq}\RRR_{pjqk}-2g^{pq}\RRR_{jp}\RRR_{kq})g_{il}\\
&\,+4\RRR_{jk}\RRR_{il}-4\RRR_{ik}\RRR_{jl}\Bigr)\\
=&\,\Delta\BBB_{ijkl}-\frac{2}{n-2}\Bigl((\RRR^{pq}\RRR_{piqk}-g^{pq}\RRR_{ip}\RRR_{kq})g_{jl}
-(\RRR^{pq}\RRR_{piql}-g^{pq}\RRR_{ip}\RRR_{lq})g_{jk}\\
&\,+(\RRR^{pq}\RRR_{pjql}-g^{pq}\RRR_{jp}\RRR_{lq})g_{ik}-(\RRR^{pq}\RRR_{pjqk}
-g^{pq}\RRR_{jp}\RRR_{kq})g_{il}\Bigr)\\
&\,+\frac{4}{n-2}(\RRR_{ik}\RRR_{jl}-\RRR_{jk}\RRR_{il})\,.
\end{align*}
Now we deal with the terms like $\RRR^{pq}\RRR_{piqk}$.\\
We have by definition
$\RRR^{pq}\RRR_{piqk}=\RRR^{pq}\WWW_{piqk}-\RRR^{pq}\AAA_{piqk}-\RRR^{pq}\BBB_{piqk}$
and
\begin{align*}
\RRR^{pq}\AAA_{piqk}=&\,\frac{\RRR}{(n-1)(n-2)}(\RRR^{pq}g_{pq}g_{ik}-\RRR^{pq}g_{pk}g_{iq})\\
=&\,\frac{\RRR}{(n-1)(n-2)}(\RRR g_{ik}-\RRR_{ik})\,,\\
\RRR^{pq}\BBB_{piqk}=&\,-\frac{1}{n-2}(\RRR^{pq}\RRR_{pq}g_{ik}-\RRR^{pq}\RRR_{pk}g_{iq}
+\RRR^{pq}\RRR_{ik}g_{pq}-\RRR^{pq}\RRR_{iq}g_{pk})\\
=&\,-\frac{1}{n-2}(\vert\Ric\vert^2 g_{ik}+\RRR\RRR_{ik}-2g^{pq}\RRR_{ip}\RRR_{kq})\,,
\end{align*}
hence, we get
\begin{align*}
\RRR^{pq}\RRR_{piqk}=&\,\RRR^{pq}\WWW_{piqk}-\frac{\RRR}{(n-1)(n-2)}(\RRR g_{ik}-\RRR_{ik})\\
&\,+\frac{1}{n-2}(\vert\Ric\vert^2
g_{ik}+\RRR\RRR_{ik}-2g^{pq}\RRR_{ip}\RRR_{kq})\\
=&\,\RRR^{pq}\WWW_{piqk}
+\frac{1}{n-2}(\vert\Ric\vert^2
g_{ik}-2g^{pq}\RRR_{ip}\RRR_{kq})
+\frac{\RRR}{(n-1)(n-2)}(n\RRR_{ik}-\RRR g_{ik})\,.
\end{align*}
Substituting these terms in the formula for $\dt \BBB_{ijkl}$ we obtain
\begin{align*}
\dt\BBB_{ijkl}
=&\,\Delta\BBB_{ijkl}
-\frac{2}{n-2}(
\RRR^{pq}\WWW_{piqk}g_{jl}
-\RRR^{pq}\WWW_{piql}g_{jk}
+\RRR^{pq}\WWW_{pjql}g_{ik}
-\RRR^{pq}\WWW_{pjqk}g_{il})\\
&\,-\frac{2\vert\Ric\vert^2}{(n-2)^2}(
g_{ik}g_{jl}
-g_{il}g_{jk}
+g_{jl}g_{ik}
-g_{jk}g_{il})\\
&\,+\frac{4}{(n-2)^2}(
g^{pq}\RRR_{ip}\RRR_{kq}g_{jl}
-g^{pq}\RRR_{ip}\RRR_{lq}g_{jk}
+g^{pq}\RRR_{jp}\RRR_{lq}g_{ik}
-g^{pq}\RRR_{jp}\RRR_{kq}g_{il})\\
&\,-\frac{2n\RRR}{(n-1)(n-2)^2}(
\RRR_{ik}g_{jl}
-\RRR_{il}g_{jk}
+\RRR_{jl}g_{ik}
-\RRR_{jk}g_{il})\\
&\,+\frac{2\RRR^2}{(n-1)(n-2)^2}(
g_{ik}g_{jl}
-g_{il}g_{jk}
+g_{jl}g_{ik}
-g_{jk}g_{il})\\
&\,+\frac{2}{n-2}(g^{pq}\RRR_{ip}\RRR_{kq}g_{jl}
-g^{pq}\RRR_{ip}\RRR_{lq}g_{jk}
+g^{pq}\RRR_{jp}\RRR_{lq}g_{ik}
-g^{pq}\RRR_{jp}\RRR_{kq}g_{il})\\
&\,+\frac{4}{n-2}(\RRR_{ik}\RRR_{jl}-\RRR_{jk}\RRR_{il})\\
=&\,\Delta\BBB_{ijkl}
-\frac{2}{n-2}(
\RRR^{pq}\WWW_{piqk}g_{jl}
-\RRR^{pq}\WWW_{piql}g_{jk}
+\RRR^{pq}\WWW_{pjql}g_{ik}
-\RRR^{pq}\WWW_{pjqk}g_{il})\\
&\,+\frac{2n}{(n-2)^2}(
g^{pq}\RRR_{ip}\RRR_{kq}g_{jl}
-g^{pq}\RRR_{ip}\RRR_{lq}g_{jk}
+g^{pq}\RRR_{jp}\RRR_{lq}g_{ik}
-g^{pq}\RRR_{jp}\RRR_{kq}g_{il})\\
&\,-\frac{2n\RRR}{(n-1)(n-2)^2}(
\RRR_{ik}g_{jl}
-\RRR_{il}g_{jk}
+\RRR_{jl}g_{ik}
-\RRR_{jk}g_{il})\\
&\,+\frac{2\RRR^2-2(n-1)\vert\Ric\vert^2}{(n-1)(n-2)^2}(
g_{ik}g_{jl}
-g_{il}g_{jk}
+g_{jl}g_{ik}
-g_{jk}g_{il})\\
&\,+\frac{4}{n-2}(\RRR_{ik}\RRR_{jl}-\RRR_{jk}\RRR_{il})\\
=&\,\Delta\BBB_{ijkl}
-\frac{2}{n-2}(
\RRR^{pq}\WWW_{piqk}g_{jl}
-\RRR^{pq}\WWW_{piql}g_{jk}
+\RRR^{pq}\WWW_{pjql}g_{ik}
-\RRR^{pq}\WWW_{pjqk}g_{il})\\
&\,+\frac{2n}{(n-2)^2}(
g^{pq}\RRR_{ip}\RRR_{kq}g_{jl}
-g^{pq}\RRR_{ip}\RRR_{lq}g_{jk}
+g^{pq}\RRR_{jp}\RRR_{lq}g_{ik}
-g^{pq}\RRR_{jp}\RRR_{kq}g_{il})\\
&\,+\frac{2n\RRR}{(n-1)(n-2)}\BBB_{ijkl}
+\frac{4\RRR}{n-2}\AAA_{ijkl}
-\frac{4\vert\Ric\vert^2}{(n-2)^2}(
g_{ik}g_{jl}-g_{il}g_{jk})\\
&\,+\frac{4}{n-2}(\RRR_{ik}\RRR_{jl}-\RRR_{jk}\RRR_{il})\,.
\end{align*}
Hence, 
\begin{align}
\Bigl(\dt-\Delta\Bigr)\WWW_{ijkl}=&\,\label{Weq}
\Bigl(\dt-\Delta\Bigr)(\RRR_{ijkl}+\AAA_{ijkl}+\BBB_{ijkl})\\
=&\,2(\CCC_{ijkl}-\CCC_{ijlk}+\CCC_{ikjl}-\CCC_{iljk})\nonumber\\
&\,-g^{pq}(\RRR_{ip}\RRR_{qjkl}+\RRR_{jp}\RRR_{iqkl}+\RRR_{kp}\RRR_{ijql}+\RRR_{lp}\RRR_{ijkq})\nonumber\\
&\,+\frac{2\vert\Ric\vert^2}{(n-1)(n-2)}(g_{ik}g_{jl}-g_{il}g_{jk})
+\frac{2\RRR}{n-1}\BBB_{ijkl}\nonumber\\
&\,-\frac{2}{n-2}(
\RRR^{pq}\WWW_{piqk}g_{jl}
-\RRR^{pq}\WWW_{piql}g_{jk}
+\RRR^{pq}\WWW_{pjql}g_{ik}
-\RRR^{pq}\WWW_{pjqk}g_{il})\nonumber\\
&\,+\frac{2n}{(n-2)^2}(
g^{pq}\RRR_{ip}\RRR_{kq}g_{jl}
-g^{pq}\RRR_{ip}\RRR_{lq}g_{jk}
+g^{pq}\RRR_{jp}\RRR_{lq}g_{ik}
-g^{pq}\RRR_{jp}\RRR_{kq}g_{il})\nonumber\\
&\,+\frac{2n\RRR}{(n-1)(n-2)}\BBB_{ijkl}
+\frac{4\RRR}{n-2}\AAA_{ijkl}
-\frac{4\vert\Ric\vert^2}{(n-2)^2}(
g_{ik}g_{jl}-g_{il}g_{jk})\nonumber\\
&\,+\frac{4}{n-2}(\RRR_{ik}\RRR_{jl}-\RRR_{jk}\RRR_{il})\nonumber\\
=&\,2(\CCC_{ijkl}-\CCC_{ijlk}+\CCC_{ikjl}-\CCC_{iljk})\nonumber\\
&\,-g^{pq}(\RRR_{ip}\RRR_{qjkl}+\RRR_{jp}\RRR_{iqkl}+\RRR_{kp}\RRR_{ijql}+\RRR_{lp}\RRR_{ijkq})\nonumber\\
&\,-\frac{2}{n-2}(
\RRR^{pq}\WWW_{piqk}g_{jl}
-\RRR^{pq}\WWW_{piql}g_{jk}
+\RRR^{pq}\WWW_{pjql}g_{ik}
-\RRR^{pq}\WWW_{pjqk}g_{il})\nonumber\\
&\,+\frac{2n}{(n-2)^2}(
g^{pq}\RRR_{ip}\RRR_{kq}g_{jl}
-g^{pq}\RRR_{ip}\RRR_{lq}g_{jk}
+g^{pq}\RRR_{jp}\RRR_{lq}g_{ik}
-g^{pq}\RRR_{jp}\RRR_{kq}g_{il})\nonumber\\
&\,-\frac{4\RRR}{(n-2)^2}(\RRR_{ik}g_{jl}-\RRR_{il}g_{jk}
+\RRR_{jl}g_{ik}-\RRR_{jk}g_{il})\nonumber\\
&\,+\frac{4\RRR^2-2n\vert\Ric\vert^2}{(n-1)(n-2)^2}(g_{ik}g_{jl}-g_{il}g_{jk})
+\frac{4}{n-2}(\RRR_{ik}\RRR_{jl}-\RRR_{jk}\RRR_{il})\,.\nonumber
\end{align}
Now, in order to simplify the formulas, we assume to be in an
orthonormal basis, then $\CCC_{ijkl}=\RRR_{pijq}\RRR_{qlkp}$ and we
have
\begin{align*}
\CCC_{ijkl}=&\,\RRR_{pijq}\RRR_{qlkp}\\
=&\,\WWW_{pijq}\WWW_{qlkp}+\AAA_{pijq}\AAA_{qlkp}+\BBB_{pijq}\BBB_{qlkp}
+\AAA_{pijq}\BBB_{qlkp}+\BBB_{pijq}\AAA_{qlkp}\\
&\,-\WWW_{pijq}\AAA_{qlkp}-\WWW_{pijq}\BBB_{qlkp}-\AAA_{pijq}\WWW_{qlkp}
-\BBB_{pijq}\WWW_{qlkp}\,.
\end{align*}
Substituting the expressions for the tensors $\AAA$ and $\BBB$ in 
the above terms and simplifying, we obtain the following identities.
$$
\AAA_{pijq}\AAA_{qlkp}=\frac{\RRR^{2}}{(n-1)^{2}(n-2)^{2}}(g_{ik}g_{jl}+(n-2)g_{ij}g_{lk})\,,
$$
\begin{align*}
\BBB_{pijq}\BBB_{qlkp}=&\,\frac{1}{(n-2)^{2}}(\RRR_{pj}g_{iq}+\RRR_{iq}g_{pj}-\RRR_{pq}g_{ij}-\RRR_{ij}g_{pq})(\RRR_{qk}g_{lp}+\RRR_{lp}g_{qk}-\RRR_{pq}g_{lk}-\RRR_{lk}g_{pq})\\
=&\,\frac{1}{(n-2)^{2}}\big(2\RRR_{ik}\RRR_{lj}+(n-4)\RRR_{ij}\RRR_{lk}+\RRR_{pj}\RRR_{pl}g_{ik}+\RRR_{pk}\RRR_{pi}g_{lj}-2\RRR_{pj}\RRR_{pi}g_{lk}-2\RRR_{pl}\RRR_{pk}g_{ij}\\
&\,+\RRR \RRR_{ij} g_{lk} + \RRR \RRR_{lk}g_{ij} +|\Ric|^{2}g_{ij}g_{lk}\big)\,,
\end{align*}
$$
\AAA_{pijq}\BBB_{qlkp}=-\frac{\RRR}{(n-1)(n-2)^{2}}\big(\RRR_{ik}g_{lj}+\RRR_{lj}g_{ik}-\RRR_{ij}g_{lk}+(n-3)\RRR_{lk}g_{ij}
+ \RRR g_{ij}g_{lk}\big)\,,
$$
$$
\BBB_{pijq}\AAA_{qlkp}=-\frac{\RRR}{(n-1)(n-2)^{2}}\big(\RRR_{lj}g_{ik}+\RRR_{ik}g_{lj}-\RRR_{lk}g_{ij}+(n-3)\RRR_{ij}g_{lk}
+ \RRR g_{ij}g_{lk}\big)\,,
$$
$$
\WWW_{pijq}\AAA_{qlkp}=\frac{\RRR}{(n-1)(n-2)}\WWW_{lijk}\,,
$$
$$
\AAA_{pijq}\WWW_{qlkp}=\frac{\RRR}{(n-1)(n-2)}\WWW_{ilkj}\,,
$$
$$
\WWW_{pijq}\BBB_{qlkp}=-\frac{1}{n-2}(\WWW_{lijp}\RRR_{pk}+\WWW_{pijk}\RRR_{lp}-\WWW_{pijq}\RRR_{pq}g_{lk})\,,
$$
$$
\BBB_{pijq}\WWW_{qlkp}=-\frac{1}{n-2}(\WWW_{ilkp}\RRR_{pj}+\WWW_{plkj}\RRR_{pi}-\WWW_{qlkp}\RRR_{pq}g_{ij})
$$
where in these last four computations we used the fact that every trace of
the Weyl tensor is null.

Interchanging the indexes and summing we get
\begin{align*}
\AAA_{pijq}\AAA_{qlkp}
-\AAA_{pijq}\AAA_{qklp}
+\AAA_{pikq}\AAA_{qljp}&\,
-\AAA_{pilq}\AAA_{qkjp}\\
=\frac{\RRR^{2}}{(n-1)^{2}(n-2)^{2}}
&\,\Bigl(g_{ik}g_{jl}+(n-2)g_{ij}g_{lk}
-g_{il}g_{jk}-(n-2)g_{ij}g_{lk}\\
&\,+g_{ij}g_{kl}+(n-2)g_{ik}g_{lj}
-g_{ij}g_{kl}-(n-2)g_{il}g_{jk}\Bigr)\\
=\frac{\RRR^{2}}{(n-1)(n-2)^{2}}
&\,(g_{ik}g_{jl}-g_{il}g_{jk})\,,
\end {align*}
\begin{align*}
\BBB_{pijq}\BBB_{qlkp}
-\BBB_{pijq}\BBB_{qklp}&
+\BBB_{pikq}\BBB_{qljp}
-\BBB_{pilq}\BBB_{qkjp}\\
=\frac{1}{(n-2)^{2}}
&\,\Bigl(2\RRR_{ik}\RRR_{lj}+(n-4)\RRR_{ij}\RRR_{lk}+\RRR_{pj}\RRR_{pl}g_{ik}+\RRR_{pk}\RRR_{pi}g_{lj}\\
&\,\phantom{\Bigl(}-2\RRR_{pj}\RRR_{pi}g_{lk}-2\RRR_{pl}\RRR_{pk}g_{ij}+\RRR \RRR_{ij} g_{lk} + \RRR \RRR_{lk}g_{ij}
+|\Ric|^{2}g_{ij}g_{lk}\\
&\,\phantom{\Bigl(}-2\RRR_{il}\RRR_{kj}-(n-4)\RRR_{ij}\RRR_{lk}-\RRR_{pj}\RRR_{pk}g_{il}-\RRR_{pl}\RRR_{pi}g_{kj}\\
&\,\phantom{\Bigl(}+2\RRR_{pj}\RRR_{pi}g_{lk}+2\RRR_{pk}\RRR_{pl}g_{ij}-\RRR \RRR_{ij} g_{lk}-\RRR \RRR_{lk}g_{ij}
-|\Ric|^{2}g_{ij}g_{lk}\\
&\,\phantom{\Bigl(}+2\RRR_{ij}\RRR_{lk}+(n-4)\RRR_{ik}\RRR_{lj}+\RRR_{pk}\RRR_{pl}g_{ij}+\RRR_{pj}\RRR_{pi}g_{lk}\\
&\,\phantom{\Bigl(}-2\RRR_{pk}\RRR_{pi}g_{lj}-2\RRR_{pl}\RRR_{pj}g_{ik}+\RRR\RRR_{ik}g_{lj}+\RRR\RRR_{lj}g_{ik}
+|\Ric|^{2}g_{ik}g_{lj}\\
&\,\phantom{\Bigl(}-2\RRR_{ij}\RRR_{kl}-(n-4)\RRR_{il}\RRR_{jk}-\RRR_{pl}\RRR_{pk}g_{ij}-\RRR_{pj}\RRR_{pi}g_{kl}\\
&\,\phantom{\Bigl(}+2\RRR_{pl}\RRR_{pi}g_{jk}+2\RRR_{pk}\RRR_{pj}g_{il}-\RRR\RRR_{il}g_{jk}-\RRR\RRR_{jk}g_{il}
-|\Ric|^{2}g_{il}g_{jk}\Bigr)\\
=\frac{1}{(n-2)^{2}}
&\,\Bigl(
(n-2)(\RRR_{ik}\RRR_{lj}-\RRR_{il}\RRR_{jk})\\
&\,\phantom{\Bigl(}-\RRR_{pj}\RRR_{pl}g_{ik}-\RRR_{pk}\RRR_{pi}g_{lj}
+\RRR_{pl}\RRR_{pi}g_{jk}+\RRR_{pk}\RRR_{pj}g_{il}\\
&\,\phantom{\Bigl(}+\RRR(\RRR_{ik}g_{lj}+\RRR_{lj}g_{ik}-\RRR_{il}g_{jk}-\RRR_{jk}g_{il})\\
&\,\phantom{\Bigl(}+|\Ric|^{2}(g_{ik}g_{lj}-g_{il}g_{jk})\Bigr)\,,
\end{align*}
\begin{align*}
\AAA_{pijq}\BBB_{qlkp}+\BBB_{pijq}\AAA_{qlkp}-&\AAA_{pijq}\BBB_{qklp}-\BBB_{pijq}\AAA_{qklp}\\
&\,+\AAA_{pikq}\BBB_{qljp}+\BBB_{pikq}\AAA_{qljp}-\AAA_{pilq}\BBB_{qkjp}-\BBB_{pilq}\AAA_{qkjp}\\
=-\frac{\RRR}{(n-1)(n-2)^{2}}
&\,\Bigl(\RRR_{ik}g_{lj}+\RRR_{lj}g_{ik}-\RRR_{ij}g_{lk}+(n-3)\RRR_{lk}g_{ij}
+ \RRR g_{ij}g_{lk}\\
&\,+ \RRR_{lj}g_{ik}+\RRR_{ik}g_{lj}-\RRR_{lk}g_{ij}+(n-3)\RRR_{ij}g_{lk}
+ \RRR g_{lk}g_{ij}\\
&\,-\RRR_{il}g_{kj}-\RRR_{jk}g_{il}+\RRR_{ij}g_{kl}-(n-3)\RRR_{kl}g_{ij}
- \RRR g_{ij}g_{kl}\\
&\,- \RRR_{kj}g_{il}-\RRR_{il}g_{kj}+\RRR_{kl}g_{ij}-(n-3)\RRR_{ij}g_{kl}
- \RRR g_{kl}g_{ij}\\
&\,+\RRR_{ij}g_{lk}+\RRR_{lk}g_{ij}-\RRR_{ik}g_{lj}+(n-3)\RRR_{lj}g_{ik}
+ \RRR g_{ik}g_{lj}\\
&\,+ \RRR_{lk}g_{ij}+\RRR_{ij}g_{lk}-\RRR_{lj}g_{ik}+(n-3)\RRR_{ik}g_{lj}
+ \RRR g_{lj}g_{ik}\\
&\,-\RRR_{ij}g_{kl}-\RRR_{lk}g_{ij}+\RRR_{il}g_{kj}-(n-3)\RRR_{kj}g_{il}
- \RRR g_{il}g_{kj}\\
&\,- \RRR_{kl}g_{ij}-\RRR_{ij}g_{kl}+\RRR_{kj}g_{il}-(n-3)\RRR_{il}g_{kj}
- \RRR g_{kj}g_{il}\Bigl)\\
=-\frac{\RRR}{(n-1)(n-2)\phantom{^{2}}}&\,\Bigl(\RRR_{ik}g_{jl}+\RRR_{jl}g_{ik}-\RRR_{jk}g_{il}-\RRR_{il}g_{jk}\Bigl)\\
-\frac{2\RRR^{2}}{(n-1)(n-2)^{2}}&\,(g_{ik}g_{jl}-g_{il}g_{jk})\,
\end {align*}
and
\begin{align*}
\WWW_{pijq}\AAA_{qlkp}&
-\WWW_{pijq}\AAA_{qklp}
+\WWW_{pikq}\AAA_{qljp}
-\WWW_{pilq}\AAA_{qkjp}\\
=&\,\frac{\RRR}{(n-1)(n-2)}(
\WWW_{lijk}-\WWW_{kijl}+\WWW_{likj}-\WWW_{kilj})=0\,,
\end{align*}
since the Weyl tensor, sharing the same symmetries of the Riemann
tensor, is skew--symmetric in the third--fourth
indexes.\\
The same result holds for the other sum as
$$
\AAA_{pijq}\WWW_{qlkp}=\frac{\RRR}{(n-1)(n-2)}\WWW_{ilkj}
=\frac{\RRR}{(n-1)(n-2)}\WWW_{lijk}=\WWW_{pijq}\AAA_{qlkp}
$$
hence, 
$$
\AAA_{pijq}\WWW_{qlkp}
-\AAA_{pijq}\WWW_{qklp}
+\AAA_{pikq}\WWW_{qljp}
-\AAA_{pilq}\WWW_{qkjp}=0\,.
$$
Finally, for the remaining two terms we have
\begin{align*}
-\WWW_{pijq}\BBB_{qlkp}-\BBB_{pijq}\WWW_{qlkp}&\,+\WWW_{pijq}\BBB_{qklp}+\BBB_{pijq}\WWW_{qklp}\\
&\,-\WWW_{pikq}\BBB_{qljp}-\BBB_{pikq}\WWW_{qljp}+\WWW_{pilq}\BBB_{qkjp}+\BBB_{pilq}\WWW_{qkjp}\\
=\frac{1}{n-2}&\,\Bigl(\WWW_{lijp}\RRR_{pk}+\WWW_{pijk}\RRR_{lp}-\WWW_{pijq}\RRR_{pq}g_{lk}\\
&\,+\WWW_{ilkp}\RRR_{pj}+\WWW_{plkj}\RRR_{pi}-\WWW_{qlkp}\RRR_{pq}g_{ij}\\
&\,-\WWW_{kijp}\RRR_{pl}-\WWW_{pijl}\RRR_{kp}+\WWW_{pijq}\RRR_{pq}g_{kl}\\
&\,-\WWW_{iklp}\RRR_{pj}-\WWW_{pklj}\RRR_{pi}+\WWW_{qklp}\RRR_{pq}g_{ij}\\
&\,+\WWW_{likp}\RRR_{pj}+\WWW_{pikj}\RRR_{lp}-\WWW_{pikq}\RRR_{pq}g_{jl}\\
&\,+\WWW_{iljp}\RRR_{pk}+\WWW_{pljk}\RRR_{pi}-\WWW_{qljp}\RRR_{pq}g_{ik}\\
&\,-\WWW_{kilp}\RRR_{pj}-\WWW_{pilj}\RRR_{kp}+\WWW_{pilq}\RRR_{pq}g_{kj}\\
&\,-\WWW_{ikjp}\RRR_{pl}-\WWW_{pkjl}\RRR_{pi}+\WWW_{qkjp}\RRR_{pq}g_{il}\Bigr)\\
=\frac{1}{n-2}&\,\Bigl(\WWW_{pilq}\RRR_{pq}g_{kj}
+\WWW_{qkjp}\RRR_{pq}g_{il}-\WWW_{pikq}\RRR_{pq}g_{jl}-\WWW_{qljp}\RRR_{pq}g_{ik}\Bigr)
\end{align*}
where we used repeatedly the symmetries of the Weyl and the Ricci
tensors.\\
Hence, summing all these terms we conclude
\begin{align}\label{ccceq}
2(\CCC_{ijkl}-&\CCC_{ijlk}+\CCC_{ikjl}-\CCC_{iljk})=2(\DDD_{ijkl}-\DDD_{ijlk}
+\DDD_{ikjl}-\DDD_{iljk})\\
&\,+\frac{2\RRR^{2}}{(n-1)(n-2)^{2}}(g_{ik}g_{jl}-g_{il}g_{jk})\nonumber\\
&\,+\frac{2}{n-2}(\RRR_{ik}\RRR_{lj}-\RRR_{il}\RRR_{jk})\nonumber\\
&\,+\frac{2}{(n-2)^{2}}(-\RRR_{pj}\RRR_{pl}g_{ik}-\RRR_{pk}\RRR_{pi}g_{lj}+\RRR_{pl}\RRR_{pi}g_{jk}+\RRR_{pk}\RRR_{pj}g_{il})\nonumber\\
&\,+\frac{2\RRR}{(n-2)^{2}}(\RRR_{ik}g_{lj}+\RRR_{lj}g_{ik}-\RRR_{il}g_{jk}-\RRR_{jk}g_{il})+\frac{2|\Ric|^{2}}{(n-2)^2}(g_{ik}g_{lj}-g_{il}g_{jk})\nonumber\\
&\,-\frac{2\RRR}{(n-1)(n-2)}(\RRR_{ik}g_{jl}+\RRR_{jl}g_{ik}-\RRR_{jk}g_{il}-\RRR_{il}g_{jk})\nonumber\\
&\,-\frac{4\RRR^{2}}{(n-1)(n-2)^{2}}(g_{ik}g_{jl}-g_{il}g_{jk})\nonumber\\
&\,+\frac{2}{n-2}(\WWW_{pilq}\RRR_{pq}g_{kj}
+\WWW_{qkjp}\RRR_{pq}g_{il}-\WWW_{pikq}\RRR_{pq}g_{jl}-\WWW_{qljp}\RRR_{pq}g_{ik})\nonumber\\
=&\,2(\DDD_{ijkl}-\DDD_{ijlk}+\DDD_{ikjl}-\DDD_{iljk})\nonumber\\
&\,+\frac{2(n-1)\vert\Ric\vert^2-2\RRR^{2}}{(n-1)(n-2)^{2}}(g_{ik}g_{jl}-g_{il}g_{jk})\nonumber\\
&\,+\frac{2}{n-2}(\RRR_{ik}\RRR_{lj}-\RRR_{il}\RRR_{jk})\nonumber\\
&\,-\frac{2}{(n-2)^{2}}(\RRR_{pj}\RRR_{pl}g_{ik}+\RRR_{pk}\RRR_{pi}g_{lj}-\RRR_{pl}\RRR_{pi}g_{jk}-\RRR_{pk}\RRR_{pj}g_{il})\nonumber\\
&\,+\frac{2\RRR}{(n-1)(n-2)^2}(\RRR_{ik}g_{jl}+\RRR_{jl}g_{ik}-\RRR_{jk}g_{il}-\RRR_{il}g_{jk})\nonumber\\
&\,+\frac{2}{n-2}(\WWW_{pilq}\RRR_{pq}g_{kj}
+\WWW_{qkjp}\RRR_{pq}g_{il}-\WWW_{pikq}\RRR_{pq}g_{jl}-\WWW_{qljp}\RRR_{pq}g_{ik})\,,\nonumber
\end{align}
where $\DDD_{ijkl}=\WWW_{pijq}\WWW_{qlkp}$.

Then we deal with the following term appearing in equation~\eqref{Weq},
\begin{align*}
\RRR_{ip}\RRR_{pjkl}+\RRR_{jp}\RRR_{ipkl}+&\,\RRR_{kp}\RRR_{ijpl}+\RRR_{lp}\RRR_{ijkp}\\
=&\,\RRR_{ip}\WWW_{pjkl}+\RRR_{jp}\WWW_{ipkl}+\RRR_{kp}\WWW_{ijpl}+\RRR_{lp}\WWW_{ijkp}\\
&\,-\frac{\RRR}{(n-1)(n-2)}\Bigl(\RRR_{ip}(g_{pk}g_{jl}-g_{pl}g_{jk})
+\RRR_{jp}(g_{ik}g_{pl}-g_{il}g_{pk})\Bigr)\\
&\,-\frac{\RRR}{(n-1)(n-2)}\Bigl(\RRR_{kp}(g_{ip}g_{jl}-g_{il}g_{jp})
+\RRR_{lp}(g_{ik}g_{jp}-g_{ip}g_{jk})\Bigr)\\
&\,+\frac{1}{n-2}(\RRR_{ip}(\RRR_{pk}g_{jl}-\RRR_{pl}g_{jk}
+\RRR_{jl}g_{pk}-\RRR_{jk}g_{pl}))\\
&\,+\frac{1}{n-2}(\RRR_{jp}(\RRR_{ik}g_{pl}-\RRR_{il}g_{pk}
+\RRR_{pl}g_{ik}-\RRR_{pk}g_{il}))\\
&\,+\frac{1}{n-2}(\RRR_{kp}(\RRR_{ip}g_{jl}-\RRR_{il}g_{jp}
+\RRR_{jl}g_{ip}-\RRR_{jp}g_{il}))\\
&\,+\frac{1}{n-2}(\RRR_{lp}(\RRR_{ik}g_{jp}-\RRR_{ip}g_{jk}
+\RRR_{jp}g_{ik}-\RRR_{jk}g_{ip}))\\
=&\,\RRR_{ip}\WWW_{pjkl}+\RRR_{jp}\WWW_{ipkl}+\RRR_{kp}\WWW_{ijpl}+\RRR_{lp}\WWW_{ijkp}\\
&\,+\frac{1}{n-2}(\RRR_{ip}\RRR_{pk}g_{jl}-\RRR_{ip}\RRR_{pl}g_{jk}
+\RRR_{jl}\RRR_{ik}-\RRR_{il}\RRR_{jk})\\
&\,+\frac{1}{n-2}(\RRR_{jl}\RRR_{ik}-\RRR_{jk}\RRR_{il}
+\RRR_{jp}\RRR_{pl}g_{ik}-\RRR_{jp}\RRR_{pk}g_{il})\\
&\,+\frac{1}{n-2}(\RRR_{kp}\RRR_{ip}g_{jl}-\RRR_{jk}\RRR_{il}
+\RRR_{ik}\RRR_{jl}-\RRR_{kp}\RRR_{jp}g_{il})\\
&\,+\frac{1}{n-2}(\RRR_{jl}\RRR_{ik}-\RRR_{lp}\RRR_{ip}g_{jk}
+\RRR_{lp}\RRR_{jp}g_{ik}-\RRR_{il}\RRR_{jk})\\
&\,-\frac{2\RRR}{(n-1)(n-2)}(
\RRR_{ik}g_{jl}
-\RRR_{il}g_{jk}
+\RRR_{jl}g_{ik}
-\RRR_{jk}g_{il})\\
=&\,\RRR_{ip}\WWW_{pjkl}+\RRR_{jp}\WWW_{ipkl}+\RRR_{kp}\WWW_{ijpl}+\RRR_{lp}\WWW_{ijkp}\\
&\,+\frac{2}{n-2}(\RRR_{ip}\RRR_{kp}g_{jl}-\RRR_{ip}\RRR_{lp}g_{jk}
+\RRR_{jp}\RRR_{lp}g_{ik}-\RRR_{jp}\RRR_{kp}g_{il})\\
&\,+\frac{4}{n-2}(\RRR_{ik}\RRR_{jl}-\RRR_{jk}\RRR_{il})\\
&\,-\frac{2\RRR}{(n-1)(n-2)}
(\RRR_{ik}g_{jl}-\RRR_{il}g_{jk}+\RRR_{jl}g_{ik}-\RRR_{jk}g_{il})\,.
\end{align*}
Inserting expression~\eqref{ccceq} and this last quantity in
equation~\eqref{Weq} we obtain
\begin{align*}
\Bigl(\dt-\Delta\Bigr)\WWW_{ijkl}
=&\,2(\DDD_{ijkl}-\DDD_{ijlk}+\DDD_{ikjl}-\DDD_{iljk})\\
&\,+\frac{2(n-1)\vert\Ric\vert^2-2\RRR^{2}}{(n-1)(n-2)^{2}}(g_{ik}g_{jl}-g_{il}g_{jk})\\
&\,+\frac{2}{n-2}(\RRR_{ik}\RRR_{lj}-\RRR_{il}\RRR_{jk})\\
&\,-\frac{2}{(n-2)^{2}}(\RRR_{pj}\RRR_{pl}g_{ik}+\RRR_{pk}\RRR_{pi}g_{lj}-\RRR_{pl}\RRR_{pi}g_{jk}-\RRR_{pk}\RRR_{pj}g_{il})\\
&\,+\frac{2\RRR}{(n-1)(n-2)^2}(\RRR_{ik}g_{jl}+\RRR_{jl}g_{ik}-\RRR_{jk}g_{il}-\RRR_{il}g_{jk})\\
&\,+\frac{2}{n-2}(\WWW_{pilq}\RRR_{pq}g_{kj}
+\WWW_{qkjp}\RRR_{pq}g_{il}-\WWW_{pikq}\RRR_{pq}g_{jl}-\WWW_{qljp}\RRR_{pq}g_{ik})\\
&\,-\RRR_{ip}\WWW_{pjkl}-\RRR_{jp}\WWW_{ipkl}-\RRR_{kp}\WWW_{ijpl}-\RRR_{lp}\WWW_{ijkp}\\
&\,-\frac{2}{n-2}(\RRR_{ip}\RRR_{kp}g_{jl}-\RRR_{ip}\RRR_{lp}g_{jk}
+\RRR_{jp}\RRR_{lp}g_{ik}-\RRR_{jp}\RRR_{kp}g_{il})\\
&\,-\frac{4}{n-2}(\RRR_{ik}\RRR_{jl}-\RRR_{jk}\RRR_{il})\\
&\,+\frac{2\RRR}{(n-1)(n-2)}
(\RRR_{ik}g_{jl}-\RRR_{il}g_{jk}+\RRR_{jl}g_{ik}-\RRR_{jk}g_{il})\\
&\,-\frac{2}{n-2}(
\RRR_{pq}\WWW_{piqk}g_{jl}
-\RRR_{pq}\WWW_{piql}g_{jk}
+\RRR_{pq}\WWW_{pjql}g_{ik}
-\RRR_{pq}\WWW_{pjqk}g_{il})\\
&\,+\frac{2n}{(n-2)^2}(
\RRR_{ip}\RRR_{kp}g_{jl}
-\RRR_{ip}\RRR_{lp}g_{jk}
+\RRR_{jp}\RRR_{lp}g_{ik}
-\RRR_{jp}\RRR_{kp}g_{il})\\
&\,-\frac{4\RRR}{(n-2)^2}(\RRR_{ik}g_{jl}-\RRR_{il}g_{jk}
+\RRR_{jl}g_{ik}-\RRR_{jk}g_{il})\\
&\,+\frac{4\RRR^2-2n\vert\Ric\vert^2}{(n-1)(n-2)^2}(g_{ik}g_{jl}-g_{il}g_{jk})
+\frac{4}{n-2}(\RRR_{ik}\RRR_{jl}-\RRR_{jk}\RRR_{il})\\
=&\,2(\DDD_{ijkl}-\DDD_{ijlk}+\DDD_{ikjl}-\DDD_{iljk})\\
&\,-(\RRR_{ip}\WWW_{pjkl}+\RRR_{jp}\WWW_{ipkl}+\RRR_{kp}\WWW_{ijpl}+\RRR_{lp}\WWW_{ijkp})\\
&\,+\frac{2(\RRR^2-\vert\Ric\vert^2)}{(n-1)(n-2)^{2}}(g_{ik}g_{jl}-g_{il}g_{jk})\\
&\,+\frac{2}{n-2}(\RRR_{ik}\RRR_{lj}-\RRR_{il}\RRR_{jk})\\
&\,+\frac{2}{(n-2)^{2}}(\RRR_{pj}\RRR_{pl}g_{ik}+\RRR_{pk}\RRR_{pi}g_{lj}-\RRR_{pl}\RRR_{pi}g_{jk}-\RRR_{pk}\RRR_{pj}g_{il})\\
&\,-\frac{2\RRR}{(n-2)^2}(\RRR_{ik}g_{jl}+\RRR_{jl}g_{ik}-\RRR_{jk}g_{il}-\RRR_{il}g_{jk})\,.
\end{align*}
Hence, we resume this long computation in the following
proposition, getting back to a standard coordinate basis.

\begin{prop}\label{wev} During the Ricci flow of an $n$--dimensional
  Riemannian manifold $(M^n,g)$, the Weyl tensor satisfies the
  following evolution equation
\begin{align*}
\Bigl(\dt-\Delta\Bigr)\WWW_{ijkl}=&\,\,\,2\,(\DDD_{ijkl}-\DDD_{ijlk}+\DDD_{ikjl}-\DDD_{iljk})\\
&\,-g^{pq}(\RRR_{ip}\WWW_{qjkl}+\RRR_{jp}\WWW_{iqkl}+\RRR_{kp}\WWW_{ijql}+\RRR_{lp}\WWW_{ijkq})\\
&\,+\frac{2}{(n-2)^2}g^{pq}(\RRR_{ip}\RRR_{qk}g_{jl}-\RRR_{ip}\RRR_{ql}g_{jk}+\RRR_{jp}
\RRR_{ql}g_{ik}-\RRR_{jp}\RRR_{qk}g_{il})\\
&\,-\frac{2\RRR}{(n-2)^2}(\RRR_{ik}g_{jl}-\RRR_{il}g_{jk}+\RRR_{jl}g_{ik}
-\RRR_{jk}g_{il})\\
&\,+\frac{2}{n-2}(\RRR_{ik}\RRR_{jl}-\RRR_{jk}\RRR_{il})
+\frac{2(\RRR^2-\vert\Ric\vert^2)}{(n-1)(n-2)^2}(g_{ik}g_{jl}-g_{il}g_{jk})\,,
\end{align*}
where $\DDD_{ijkl}=g^{pq}g^{rs}\WWW_{pijr}\WWW_{slkq}$.
\end{prop}

From this formula we immediately get the following rigidity result on
the eigenvalues of the Ricci tensor.

\begin{cor}\label{princ} Suppose that under the Ricci flow of $(M^n,g)$ of
  dimension $n\geq 4$, the Weyl tensor
  remains identically zero. Then, at every point, either the Ricci
  tensor is proportional to the metric or it has an eigenvalue of
  multiplicity $(n-1)$ and another of multiplicity 1.
\end{cor}
\begin{proof}
By the above proposition, as every term containing the Weyl tensor is
zero, the following relation holds at every point in space and time
\begin{align*}
0=&\,\frac{2}{(n-2)^2}g^{pq}(\RRR_{ip}\RRR_{qk}g_{jl}-\RRR_{ip}\RRR_{ql}g_{jk}
+\RRR_{jp}\RRR_{ql}g_{ik}-\RRR_{jp}\RRR_{qk}g_{il})\\
&\,+\frac{2\RRR^2}{(n-1)(n-2)^2}(g_{ik}g_{jl}-g_{il}g_{jk})
-\frac{2\vert\Ric\vert^2}{(n-1)(n-2)^2}(g_{ik}g_{jl}-g_{il}g_{jk})\\
&\,-\frac{2\RRR}{(n-2)^2}(\RRR_{ik}g_{jl}-\RRR_{il}g_{jk}
+\RRR_{jl}g_{ik}-\RRR_{jk}g_{il})+\frac{2}{n-2}(\RRR_{ik}\RRR_{jl}-\RRR_{jk}\RRR_{il})\,.
\end{align*}
In normal coordinates such that the Ricci tensor is diagonal we get,
for every couple of different eigenvectors $v_i$ with relative eigenvalues $\lambda_i$,
\begin{equation}\label{auto}
(n - 1)[\lambda_i^2 + \lambda_j^2 ] - (n - 1) \RRR ( \lambda_i +
\lambda_j ) + (n - 1)(n - 2) \lambda_i\lambda_j + \RRR^2 - \vert\Ric\vert^2 = 0\,.
\end{equation}
As $n\geq 4$, fixing $i$, then the equation above is a second
order polynomial in $\lambda_j$, hence it can only have at most 2
solutions, hence, we can conclude that there are at most three
possible values for the eigenvalues of the Ricci tensor.\\
Since the dimension is at least four, at least one eigenvalues must
have multiplicity two, let us say $\lambda_i$, hence the
equation~\eqref{auto} holds also for $i=j$, and it remains at most
only {\em one} possible value for the other eigenvalues $\lambda_l$
with $l\not=i$. In conclusion, either the eigenvalues are all equal or
they divide in only two possible values, $\lambda$ with multiplicity
larger than one, say $k$ and $\mu\not=\lambda$. Suppose that $\mu$ also has multiplicity
larger than one, that is, $k<n-1$, then we have
\begin{gather}
n\lambda^2 - 2\RRR\lambda=\frac{\vert\Ric\vert^2 -\RRR^2}{n-1}\label{111}\\
n\mu^2 - 2\RRR\mu=\frac{\vert\Ric\vert^2 -\RRR^2}{n-1}\nonumber
\end{gather}
taking the difference and dividing by $(\lambda-\mu)$ we get
$$
n(\lambda+\mu)=2\RRR=2[k\lambda+(n-k)\mu]
$$
then,
$$
(n-2k)\lambda=(n-2k)\mu
$$
hence, $n=2k$, but then getting back to equation~\eqref{111}, $\RRR=n(\mu+\lambda)/2$ and
$$
n\lambda^2 -
n(\mu+\lambda)\lambda=\frac{n(\lambda^2+\mu^2)/2-n^2(\mu^2+\lambda^2+2\lambda\mu)/4}{n-1}
$$
which implies
$$
-4n\lambda\mu=-\frac{n(n-2)}{n-1}(\lambda^2+\mu^2)-\frac{2n^2}{n-1}\mu\lambda
$$
that is, after some computation,
$$
\frac{2n(n-2)}{n-1}\mu\lambda=\frac{n(n-2)}{n-1}(\lambda^2+\mu^2)\,,
$$
which implies $\lambda=\mu$.

At the end we conclude that at every point of $M^n$, either $\Ric=\lambda g$ or there is an
eigenvalue $\lambda$ of multiplicity $(n-1)$ and another $\mu$ of
multiplicity 1.
\end{proof}

\begin{rem}
Notice that in dimension three equation~\eqref{auto} becomes
\begin{align*}
2[ \lambda_i^2 +\lambda_j^2 ] - &\,2\RRR ( \lambda_i +
\lambda_j ) + 2\lambda_i\lambda_j + \RRR^2 - \vert\Ric\vert^2\\ 
=&\, 2(\lambda_i+\lambda_j)^2 - 2\RRR ( \lambda_i +
\lambda_j ) -2\lambda_i\lambda_j + \RRR^2 - \vert\Ric\vert^2\\
=&\,-2\lambda_l(\lambda_i+\lambda_j) - 2\lambda_i\lambda_j + \RRR^2 -
\vert\Ric\vert^2\\
=&\,0\,,
\end{align*}
where $\lambda_i$, $\lambda_j$ and $\lambda_l$ are the three
eigenvalues of the Ricci tensor.\\ 
Hence, the condition is void and our argument does not work. This is
clearly not unexpected as the Weyl tensor is identically zero for every 
three--dimensional Riemannian manifold.
\end{rem}

\section{Locally Conformally Flat Ricci Solitons}

Let $(M^n,g)$, for $n\geq 4$, be a connected, complete, 
Ricci soliton, that is, there exists a smooth $1$--form $\omega$
and a constant $\alpha\in\mathbb{R}$ such that
\begin{equation*}
\RRR_{ij}+\frac{1}{2}(\nabla_{i}\omega_{j}+\nabla_{j}\omega_{i})=\frac{\alpha}{n} g_{ij}\,.
\end{equation*}
If $\alpha>0$ we say that the soliton is {\em shrinking}, if
$\alpha=0$ {\em steady}, if $\alpha<0$ {\em expanding}.\\
If there exists a smooth function $f:M^n\to\RR$ such that $df=\omega$ 
we say that the soliton is a {\em gradient} Ricci soliton and $f$ its
{\em potential function}, then we have
\begin{equation*}
\RRR_{ij}+\nabla^2_{ij}f=\frac{\alpha}{n} g_{ij}\,.
\end{equation*}

\smallskip

\noindent {\em If the metric dual field to the form $\omega$ is complete, then a
Ricci soliton generates a self--similar solution to the Ricci flow (if the
soliton is a {\em gradient} soliton this condition is automatically
satisfied~\cite{zhang2}).\\
In all this section we will assume to be in this case.}

\smallskip 

In this section we discuss the classification of Ricci solitons
$(M^n,g)$, for $n\geq 4$, which are locally conformally flat
(LCF). As a consequence of Corollary~\ref{princ} we have the following
fact.

\begin{prop}\label{eigen} Let $(M^n,g)$ be a complete, LCF Ricci soliton 
  of dimension $n\geq 4$. Then, at every point, either the Ricci
  tensor is proportional to the metric or it has an eigenvalue of
  multiplicity $(n-1)$ and another of multiplicity 1.
\end{prop}

If a manifold $(M^n,g)$ is LCF, it follows that
\begin{align*}
0=&\,\nabla^l\WWW_{ijkl}\\
=&\,\nabla^l\Bigl(\RRR_{ijkl}+\frac{\RRR}{(n-1)(n-2)}(g_{ik}g_{jl}-g_{il}g_{jk})
- \frac{1}{n-2}(\RRR_{ik}g_{jl}-\RRR_{il}g_{jk}
+\RRR_{jl}g_{ik}-\RRR_{jk}g_{il})\Bigr)\\
=&\,-\nabla_i\RRR_{jk}+\nabla_j\RRR_{ik}
+\frac{\nabla_j\RRR}{(n-1)(n-2)}g_{ik}
-\frac{\nabla_i\RRR}{(n-1)(n-2)}g_{jk}\\
&\,- \frac{1}{n-2}(\nabla_j\RRR_{ik}-\nabla^l\RRR_{il}g_{jk}
+\nabla^l\RRR_{jl}g_{ik}-\nabla_i\RRR_{jk}g_{il})\\
=&\,-\frac{n-3}{n-2}(\nabla_i\RRR_{jk}-\nabla_j\RRR_{ik})
+\frac{\nabla_j\RRR}{(n-1)(n-2)}g_{ik}
-\frac{\nabla_i\RRR}{(n-1)(n-2)}g_{jk}\\
&\,+\frac{1}{2(n-2)}(\nabla_i\RRR g_{jk}/2
-\nabla_j\RRR g_{ik}/2)\\
=&\,-\frac{n-3}{n-2}\Bigl[
\nabla_i\RRR_{jk}+\nabla_j\RRR_{ik}-\frac{(\nabla_i\RRR g_{jk}-\nabla_j\RRR g_{ik})}{2(n-1)}\Bigr]\\
=&\,\frac{n-3}{n-2} \Bigl[ \nabla_{j} \Big(\RRR_{ik}
  - \frac{1}{2(n-1)}\RRR g_{ik}\Big) - \nabla_{i}\Big(\RRR_{jk} -
  \frac{1}{2(n-1)}\RRR g_{jk}\Big)\Bigr]\,,
\end{align*}
where we used the second Bianchi identity and Schur's Lemma
$\nabla\RRR=2\div\Ric$.\\
Hence, since we assumed that the dimension $n$ is at least four, the
Schouten tensor defined by
$\SSS=\Ric-\frac{1}{2(n-1)}\RRR g$ satisfies the equation
$$
(\nabla_{X}\SSS)\,Y=(\nabla_{Y}\SSS)\,X, \quad X,Y\in TM\,.
$$
Any symmetric two tensor satisfying this condition is called a Codazzi
tensor (see~\cite[Chapter~16]{besse} for a general overview of Codazzi
tensors).\\
Suppose that we have a local orthonormal
frame $\{E_{1},\dots,E_{n}\}$ in an open subset $\Omega$ of $M^n$ such that $\Ric(E_{1})=\lambda E_{1}$
and $\Ric(E_{i})=\mu E_{i}$ for $i=2,\dots,n$ and $\lambda\not=\mu$. 
For every point in $\Omega$ also the Schouten tensor $\SSS$ has
two distinct eigenvalues $\sigma_{1}$ of multiplicity one and
$\sigma_{2}$ of multiplicity $(n-1)$, with the same eigenspaces of
$\lambda$ and $\mu$ respectively, and
$$
\sigma_{1}=\frac{2n-3}{2(n-1)}\lambda-\frac{1}{2}\mu\quad\text{ and }\quad
\sigma_{2}=\frac{1}{2}\mu-\frac{1}{2(n-1)}\lambda\,.
$$ 
Splitting results for Riemannian manifolds admitting a Codazzi tensor
with only two distinct eigenvalues were obtained by Derdzinski~\cite{derdz3} and
Hiepko--Reckziegel~\cite{hiepko,hiepkoreck} 
(see again~\cite[Chapter~16]{besse} for further
discussion). In particular, it can be proved that, if the two distinct
eigenvalues $\sigma_{1}$ and $\sigma_{2}$ are both ``constant along
the eigenspace $span\{E_{2},\dots,E_{n}\}$'' then the manifold is
locally a warped product on an interval of $\RR$ of a
$(n-1)$--dimensional Riemannian manifold (see~\cite[Chapter~16]{besse} and~\cite{toj1}).\\
Since $\sigma_{2}$ has multiplicity $(n-1)$, larger than 2, we have for any two distinct
indexes $i,j\geq 2$,
\begin{align*}
\partial_i\sigma_2=&\,\partial_i \SSS(E_j,E_j)\\
=&\,\nabla_i\SSS_{jj}+2\SSS(\nabla_{E_i}E_j,E_j)\\
=&\,\nabla_j\SSS_{ij}+2\sigma_2 g(\nabla_{E_i}E_j,E_j)\\
=&\,\partial_j\SSS(E_i,E_j) - \SSS(\nabla_{E_j}E_i,E_j) -
\SSS(E_i,\nabla_{E_j}E_j)\\
=&\, - \sigma_2g(\nabla_{E_j}E_i,E_j) - \sigma_2g(E_i,\nabla_{E_j}E_j)\\
=&\,0\,,
\end{align*}
hence, $\sigma_2$ is always constant along the
eigenspace $span\{E_{2},\dots,E_{n}\}$. The eigenvalue $\sigma_1$
instead, for a general LCF manifold, can vary, for example $\RR^{n}$
endowed with the metric
$$
g=\frac{dx^2}{[1+(x_1^2+x_2^2+\dots+x_{n-1}^2)]^2}
$$
is LCF and
$$
\RRR^g_{ij}=-(n-2)(\nabla^2_{ij}\log{A}-\nabla_i\log{A}\nabla_j\log{A})+(\Delta\log{A}-(n-2)\vert\nabla\log{A}\vert^2)\delta_{ij}
$$
where the derivatives are the standard ones of $\RR^n$ and
$A(x)=1+(x_1^2+x_2^2+\dots+x_{n-1}^2)$ 
(see~\cite[Theorem~1.159]{besse}). Hence, this Ricci tensor 
``factorizes'' on the eigenspaces $\langle
e_{1},\dots,e_{n-1}\rangle$ and $\langle e_n\rangle$ but the 
eigenvalue $\sigma_1$ of the Schouten tensor, which is given by
\begin{align*}
\sigma_1=g^{nn}\RRR^g_{nn}=&\,(\Delta\log{A}-(n-2)\vert\nabla\log{A}\vert^2)A^2\\
=&\,A\Delta A-(n-1)\vert\nabla A\vert^2\\
=&\,2(n-1)A-4(n-1)(A-1)\\
=&\,-2(n-1)(A-2)\,,
\end{align*}
is clearly not constant along the directions $e_{1},\dots,e_{n-1}$.

The best we can say in general is that the metric of $(M^n,g)$ locally
around every point can be written as $I\times N$ and 
$$
g(t,p)=\frac{dt^2+\sigma^K(p)}{[\alpha(t)+\beta(p)]^2}
$$
where $\sigma^K$ is a metric on $N$ of constant curvature $K$,
$\alpha:I\to\RR$ and $\beta:N\to\RR$ are smooth functions such that
${\mathrm{Hess}}^K\beta=f\sigma^K$, for some function $f:N\to\RR$ and 
where ${\mathrm{Hess}^K}$ is the Hessian of $(N,\sigma^K)$.

\subsection{Compact LCF Ricci Solitons}

A compact Ricci soliton is actually a gradient soliton (by the
work of Perelman~\cite{perel1}).\\
In general (even if they are not LCF), steady and expanding compact Ricci
solitons are Einstein, hence, when also LCF, they are of constant curvature
(respectively zero and negative).\\
In~\cite{caowang,mantemin2} it is proved that also shrinking, compact,
LCF Ricci solitons are of constant positive curvature, hence quotients
of spheres.  

\medskip

{\em Any compact, $n$--dimensional, LCF Ricci soliton is a 
quotient of $\RR^n$, $\SS^n$ and $\HH^n$ with their 
canonical metrics, for every $n\in\NN$.}

\medskip

\subsection{LCF Ricci Solitons with Constant Scalar Curvature}

Getting back to the Schouten tensor, 
if the scalar curvature $\RRR$ of an LCF Ricci soliton $(M^n,g)$ is constant, we have
that also the other eigenvalue $\sigma_1$ of the Schouten tensor is constant
along the eigenspace $span\{E_{2},\dots,E_{n}\}$, that is,
$\partial_i\sigma_1=0$, by simply differentiating the equality
$\RRR=\frac{2(n-1)}{n-2}(\sigma_1+(n-1)\sigma_2)$.

Hence, by the above discussion, we can 
conclude that around every point of $M^n$ in the open set
$\Omega\subset M^n$ where the two eigenvalues of
the Ricci tensor are distinct the manifold $(M^n,g)$ is locally a warped
product $I\times N$ with $g(t,p)=dt^2+ h^2(t)\sigma(p)$ (this
argument is due to Derdzinski~\cite{derdz3}).\\
Then the LCF hypothesis implies that
the warp factor $(N,\sigma)$ is actually a space of constant curvature
$K$ (see for instance~\cite{brovargaz}).

As the scalar curvature $\RRR$ is constant, by the evolution equation
$\partial_t\RRR=\Delta\RRR+2\vert\Ric\vert^2$ we see that also
$\vert\Ric\vert^2$ is constant, that is, locally
$\RRR=\lambda+(n-1)\mu=C_1$ and 
$\vert\Ric\vert^2=\lambda^2+(n-1)\mu^2=C_2$. Putting together these
two equations it is easy to see that then both the eigenvalues $\mu$ and
$\lambda$ are locally constant in $\Omega$.
Hence, by connectedness, either $(M^n,g)$ is Einstein, so a constant curvature space, or
the Ricci tensor has two distinct constant eigenvalues everywhere.
Using now the local warped product representation, the
Ricci tensor is expressed by (see~\cite[Proposition~9.106]{besse} or~\cite[p.~65]{chowluni}
or~\cite[p.~168]{bry2})
\begin{equation}\label{ricwarp}
\Ric=-(n-1)\,\frac{h^{\prime\prime}}{h}dt^{2}+\big((n-2)K-h\,h^{\prime\prime}-(n-2)({h^\prime})^{2}\big)\sigma^{K}\,.
\end{equation}
hence, $h^{\prime\prime}/h$ and $((n-2)K-h\,h^{\prime\prime}-(n-2)({h^\prime})^{2})/h^2$ are constant
in $t$. This implies that $(K-({h^\prime})^{2})/h^2$ is also constant
and $h^{\prime\prime}=Ch$, then locally either the manifold $(M^n,g)$ is of
constant curvature or it is the Riemannian product of a constant curvature space
with an interval of $\RR$.\\
By a maximality argument, passing to the universal
covering of the manifold, we get the following conclusion.

\medskip

{\em If $n\geq 4$, any $n$--dimensional, LCF Ricci soliton 
with constant scalar curvature is either a
quotient of $\RR^n$, $\SS^n$ and $\HH^n$ with their 
canonical metrics or a quotient of the Riemannian
products $\RR\times\SS^{n-1}$ and $\RR\times\HH^{n-1}$
(see also~\cite{pw}).}

\medskip

\subsection{Gradient LCF Ricci Solitons with Nonnegative Ricci Tensor}

Getting back again to the Codazzi property of the Schouten tensor $\SSS$,
for every index $i>1$, we have locally
$$
0=\nabla_{1}\RRR_{i1}-\nabla_{i}\RRR_{11}
  - \frac{\partial_1\RRR}{2(n-1)}g_{i1}
+\frac{\partial_i\RRR}{2(n-1)}g_{11}=
\nabla_{1}\RRR_{i1}-\nabla_{i}\RRR_{11}
+\frac{\partial_i\RRR}{2(n-1)}\,.
$$
If the soliton is a gradient LCF Ricci soliton, that is, 
$\Ric=-\nabla^2f+\frac{\alpha}{n} g$, we have $\RRR=-\Delta f+\alpha$
and taking the divergence of both sides
\begin{align*}
\partial_i\RRR/2=&\,\div\Ric_i\\
=&\,g^{jk}\nabla_k\RRR_{ij}\\
=&\,-g^{jk}\nabla_k\nabla_i\nabla_jf\\
=&\,-g^{jk}\nabla_i\nabla_k\nabla_jf-g^{jk}\RRR_{kijl}\nabla^lf\\
=&\,-\nabla_i\Delta f-\RRR_{il}\nabla^lf\\
=&\,\partial_i\RRR-\RRR_{il}\nabla^lf\,,
\end{align*}
where we used Schur's  Lemma
$\partial_i\RRR=2\div\Ric_i$ and the formula for the interchange of covariant
derivatives.\\
Hence, the relation 
$\partial_i\RRR=2\RRR_{il}\nabla^lf$ holds and 
$$
\nabla_{1}\nabla^2_{i1}f-\nabla_{i}\nabla^2_{11}f=\frac{\RRR_{ij}\nabla^jf}{n-1}\,.
$$
By means of the fact that $\WWW=0$, we compute now for $i>1$ (this is a
special case of the computation in Lemma~3.1 of~\cite{caochen}),
\begin{align*}
\frac{\mu}{n-1}\nabla_i f=&\,
\frac{\RRR_{ij}\nabla^jf}{n-1}\\
=&\,\nabla_{1}\nabla^2_{i1}f-\nabla_{i}\nabla^2_{11}f\\
=&\,\RRR_{1i1j}\nabla^{j}f\\
=&\,\left[\frac{1}{n-2}(\RRR_{11}g_{ij}-\RRR_{1j}g_{i1}
+\RRR_{ij}g_{11}-\RRR_{i1}g_{1j})
-\frac{\RRR}{(n-1)(n-2)}(g_{11}g_{ij}-g_{1j}g_{i1})\right]\nabla^jf\\
=&\,\left[\frac{1}{n-2}(\lambda g_{ij}+\mu g_{ij})
-\frac{\RRR}{(n-1)(n-2)}g_{ij}\right]\nabla^jf \\
=&\,\left[\frac{\lambda+\mu}{n-2}
-\frac{\RRR}{(n-1)(n-2)}\right]\nabla_if \\
=&\,\frac{(n-1)\lambda+(n-1)\mu -\lambda-(n-1)\mu}{(n-1)(n-2)}\nabla_if\\
=&\,\frac{\lambda}{n-1}\nabla_if\,.
\end{align*}
Then, in the open set
$\Omega\subset M^n$ where the two eigenvalues of
the Ricci tensor are distinct, the vector field $\nabla f$ is
parallel to $E_1$, hence it is an eigenvector of the Ricci tensor and
$\partial_i\RRR=2\RRR_{il}\nabla^lf=0$, for every index $i>1$.\\
As $\sigma_1=\frac{n-2}{2(n-1)}\RRR-(n-1)\sigma_2$ we get that 
also $\partial_i\sigma_1=0$ for every index $i>1$.

The set $\Omega$ is dense, otherwise its complement where
$\Ric-Rg/n=0$ has interior points and, by
Schur's Lemma, the scalar curvature would be constant in some open set
of $M^n$. Then, strong maximum principle applied to the equation 
$\partial_t\RRR=\Delta\RRR+2\vert\Ric\vert^2$ implies that $\RRR$
is constant everywhere on $M^n$, and we are in the previous case.

So we can conclude also in this case by the previous 
argument that the manifold, locally around every point in $\Omega$, is a warped
product on an interval of $\RR$ of a constant curvature
space $\mathbb{L}^K$. Moreover, $\Omega$ is obviously invariant by
``translation'' in the  $\mathbb{L}^K$--direction.

We consider a point $p\in\Omega$ and the maximal geodesic curve
$\gamma(t)$ passing from $p$ orthogonal 
to $\mathbb{L}^K$, contained in $\Omega$. It is easy to see that for every compact, connected
segment of such geodesic we have a
neighborhood $U$ and a representation of the metric in $g$ as
$$
g=dt^2+ h^2(t)\sigma^K\,,
$$
covering the segment with the local charts and possibly shrinking them
in the orthogonal directions.\\
Assuming from now on that the the Ricci tensor is
nonnegative, by the local warped representation
formula~\eqref{ricwarp} 
we see that $h^{\prime\prime}\leq 0$ along such geodesic, as
$\RRR_{tt}\geq 0$.\\
If such geodesic has no ``endpoints'', being concave the function
$h$ must be constant and we have either a flat quotient of $\RR^{n}$
or the Riemannian product of $\RR$ with a quotient of $\SS^{n-1}$. The same holds if the function 
$h$ is constant in some interval, indeed, the manifold would be locally a
  Riemannian product and the scalar curvature would be locally constant
  (hence we are in the case above).\\
If there is at least one endpoint, one of the following two
  situations happens:
\begin{itemize}
\item the function $h$ goes to zero at such endpoint, 
\item the geodesic hits the boundary of $\Omega$.
\end{itemize}
If $h$ goes to zero at an endpoint, by concavity $(h^\prime)^2$ must converge to some
positive limit and by the smoothness of the manifold, considering
again formula~\eqref{ricwarp}, 
the quantity $K-({h^\prime})^2$ must go to zero as $h$ goes to
zero, hence $K>0$ and the constant curvature
space $\mathbb{L}^K$ must be a quotient of the sphere $\SS^{n-1}$
(if the same happens also at the other endpoint, the manifold is compact). 
Then, by topological reasons we conclude that 
actually the only possibility for 
$\mathbb{L}^K$ is the sphere $\SS^{n-1}$ itself.\\
Assuming instead that $h$ does not go to zero at any endpoint, where the 
geodesic hits the boundary of $\Omega$ the Ricci tensor is proportional to the
metric, hence, again by the representation formula~\eqref{ricwarp}, the quantity
$K-({h^\prime})^2$ is going to zero and either $K=0$ or 
$K>0$.\\
The case $K=0$ is impossible, indeed $h^\prime$ would tend to zero at such endpoint, then by the 
concavity of $h$ the function $h^\prime$ has a sign, otherwise $h$ is
constant in an interval, implying that in some open set $(M^n,g)$ is flat, which cannot
happens since we are in $\Omega$. Thus, being $h^\prime\not=0$, $h$ concave and we
assumed that $h$ does not go to zero, there must be another endpoint
where the geodesic hits the boundary of $\Omega$, which is in
contradiction with $K=0$ since also in this point $K-({h^\prime})^2$
must go to zero but instead $h^\prime$ tends to some nonzero value.
Hence, $K$ must be positive and also in this case we are dealing
with a warped product of a quotient of $\SS^{n-1}$ on an interval of $\RR$.

Resuming, in the non--product situation, every connected piece of 
$\Omega$ is a warped product of a quotient of the sphere $\SS^{n-1}$ on some
intervals of $\RR$. Then, we can conclude that the universal cover
$(\widetilde{M},\widetilde{g})$ can be recovered ``gluing together'',
along constant curvature spheres, warped product pieces that can be
topological ``caps'' (when $h$ goes to
zero at an endpoint) and ``cylinders''. Nontrivial quotients $(M,g)$ of
$(\widetilde{M},\widetilde{g})$ are actually possible only when there
are no ``caps'' in this gluing procedure. In such case, by its concavity, the function $h$
must be constant along every piece of geodesic and the manifold 
$(\widetilde{M},\widetilde{g})$ is a Riemannian product. If there is at least one
``cap'', the whole manifold is a warped product of $\SS^{n-1}$ on an
interval of $\RR$.

\begin{rem}\label{remm} We do not know if the condition on $(M^n,g)$ to be a
  {\em gradient} LCF Ricci soliton is actually necessary to have locally a
  warped product. We conjecture that such conclusion should hold also for
{\em nongradient} LCF Ricci solitons.
\end{rem}

\medskip

{\em If $n\geq 4$, any $n$--dimensional, LCF gradient Ricci soliton with
  nonnegative Ricci tensor is either a quotient of $\RR^n$ and $\SS^n$ with their 
canonical metrics, or a quotient of $\RR\times\SS^{n-1}$
or it is a warped product of $\SS^{n-1}$ on a proper interval of $\RR$.}

\medskip

\subsection{The Classification of Steady and Shrinking Gradient LCF Ricci
  Solitons}

The class of solitons with nonnegative Ricci tensor is particularly
interesting as it includes all the shrinking and steady Ricci
solitons.\\
Indeed, by the same arguments of~\cite{zhang} (keeping in mind, in following the
proof of the main Proposition~3.2 in such paper, that the
nonnegativity of the scalar curvature for every complete, ancient
Ricci was proved in~\cite[Corollary~2.5]{chen2}), where the author generalizes 
the well--known Hamilton--Ivey curvature estimate to locally conformally flat, gradient,
shrinking Ricci solitons (Corollary~3.3 in the same paper~\cite{zhang}), it follows that actually {\em every} 
complete ancient solution $g(t)$ to the Ricci flow whose Weyl tensor
is identically zero for all times, is forced to have nonnegative
curvature operator for every time $t$.\\
In particular, this holds for any complete, steady or shrinking Ricci soliton (even if
not gradient) as they generate self--similar ancient solutions of Ricci flow.

By the previous discussion and the analysis
of Bryant in the steady case~\cite{bry2} (see also~\cite[Chapter~1,
Section~4]{chowbookI}) showing that there exists a unique (up to
dilation of the metric) nonflat,
steady, gradient Ricci soliton which is a warped
product of $\SS^{n-1}$ on a halfline of $\RR$, called {\em Bryant
  soliton}, we get the following classification.

\begin{prop} The steady, gradient, LCF Ricci solitons of dimension
  $n\geq4$ are given by the quotients of $\RR^n$ and the Bryant soliton.
\end{prop}

This classification result, including also the three--dimensional LCF
case, was first obtained recently by H.-D.~Cao and Q.~Chen~\cite{caochen}.

\smallskip

In the shrinking case, the analysis of Kotschwar~\cite{kotschwar} of
rotationally invariant shrinking, gradient Ricci solitons gives the
following classification where the {\em
  Gaussian soliton} is defined as the flat $\RR^n$ with a potential
function $f=\alpha\vert x\vert^2/2n$, for a constant $\alpha\in\RR$.

\begin{prop} The shrinking, gradient, LCF Ricci solitons of dimension
  $n\geq4$ are given
  by the quotients of $\SS^n$, the Gaussian solitons with $\alpha>0$
  and quotients of $\RR\times\SS^{n-1}$.
\end{prop}

This classification of shrinking, gradient, LCF Ricci solitons follows
by the works of L.~Ni and N.~Wallach~\cite{nw2}, P.~Petersen and
W.~Wylie~\cite{pw} and Z.-H.~Zhang~\cite{zhang}.\\
Several other authors contributed to the
subject, including  X.~Cao, B.~Wang and Z.~Zhang~\cite{caowang},
B.-L.~Chen~\cite{chen2}, M.~Fern\'andez--L\'opez and
E.~Garc\'{\i}a--R\'{\i}o~\cite{fergarciario2}, M.~Eminenti, G.~La~Nave and C.~Mantegazza~\cite{mantemin2}, 
O.~Munteanu and N.~Sesum~\cite{sesummunte} and again P.~Petersen and W.~Wilye~\cite{pw1}.

\smallskip

We show now that every complete, warped, LCF Ricci soliton with
nonnegative Ricci tensor is actually a gradient
soliton.\\
Proving our conjecture in Remark~\ref{remm} that 
every Ricci soliton is locally a warped product would then lead to
have a general classification of also nongradient Ricci solitons, in
the steady and shrinking cases.

\begin{rem} In the compact case, the fact that every Ricci soliton is actually a gradient
is a consequence of the work of
Perelman~\cite{perel1}. Naber~\cite{na1} showed that it is true also
for shrinking Ricci solitons with bounded curvature.\\
For examples of nongradient Ricci solitons see Baird and Danielo~\cite{bairdan}.
\end{rem}

\begin{prop}\label{grad} Let $(M^n,g)$ be a complete, warped, 
  LCF Ricci soliton with nonnegative Ricci tensor, then it is a gradient
  Ricci soliton with a potential function $f:M^n\to\RR$ depending only on the
  $t$ variable of the warping interval.
\end{prop}
\begin{proof}
We assume that $(M^n,g)$ is globally described by $M^n=I\times\mathbb{L}^K$ and 
$$
g=dt^2+ h^2(t)\sigma^K\,,
$$
where $I$ is an interval of $\RR$ or $\SS^1$ and
$(\mathbb{L}^K,\sigma^K)$ is a complete space of constant curvature $K$.\\
In the case $h$ is constant, which clearly follows if $I=\SS^1$, as
$h^{\prime\prime}\leq 0$ the conclusion is trivial.\\
We deal then with the case where $h:I\to\RR$ is zero at some
point, let us say $h(0)=0$ and $I=[0,+\infty)$, (if the interval $I$ is
bounded the manifold $M^n$ is compact and we are done). Then, $\mathbb{L}^K=\SS^{n-1}$ with its constant
curvature metric $\sigma^K$. As a consequence, we have $M^n=\RR^n$, simply connected.
We consider the form $\omega$ satisfying the structural equation
\begin{equation*}
\RRR_{\gamma\beta}+\frac{1}{2}(\nabla_\gamma\omega_\beta+\nabla_\beta\omega_\gamma)
=\frac{\alpha}{n} g_{\gamma\beta}\,,
\end{equation*}
If $\varphi:\SS^{n-1}\to\SS^{n-1}$ is an isometry of the standard
sphere, the associated map $\phi:M^n\to M^n$ given by $\phi(t,p)=(t,\varphi(p))$
is also an isometry, moreover, by the warped structure of $M^n$ we
have that the 1--form $\phi^*\omega$ also satisfies
\begin{equation*}
\RRR_{\gamma\beta}+\frac{1}{2}\big[(\nabla\phi^*\omega)_{\gamma\beta}+(\nabla\phi^*\omega)_{\beta\gamma}\big]=\frac{\alpha}{n} g_{\gamma\beta}\,,
\end{equation*}
Calling $\mathcal{I}$ the Lie group of isometries of $\SS^{n-1}$ and $\xi$ the
Haar unit measure associated to it, we 
define the following 1--form
\begin{equation*}
\theta=\int_{\mathcal{I}} \phi^*\omega\,d\xi(\varphi)\,.
\end{equation*}
By the linearity of the structural equation, we have
\begin{equation*}
\RRR_{\gamma\beta}+\frac{1}{2}(\nabla_\gamma\theta_\beta+\nabla_\beta\theta_\gamma)
=\frac{\alpha}{n} g_{\gamma\beta}\,,
\end{equation*}
moreover, by construction, we have $L_X\theta=0$ for every vector
field $X$ on $M^n$ which is a generator of an isometry $\phi$ of $M^n$
as above (in other words, $\theta$ depends only on $t$). Computing in
normal coordinates on $\SS^{n-1}$, we get
\begin{align*}
\nabla_i\theta_j=&\,-\theta(\nabla_j\partial_i)=-\Gamma_{ij}^t\theta_t=hh^\prime \sigma^K_{ij}\theta_t\,,\\
\nabla_i\theta_t=&\,-\theta(\nabla_t\partial_i)=-\Gamma_{ti}^j\theta_j=-\frac{h^\prime}{h}\theta_i\,.\\
\end{align*}
Hence,
\begin{align*}
\frac{\alpha}{n}=&\,\RRR_{tt}+\nabla_{t}\theta_{t}=-(n-1)\frac{h^{\prime\prime}}{h}
+\partial_{t}\theta_{t}\,,\\
0=&\,\nabla_{i}\theta_{t}+\nabla_{t}\theta_{i}=\partial_{t}\theta_{i}-2\frac{h^\prime}{h}\theta_i\,,\\
\frac{\alpha}{n}
g_{ij}^K=&\,\RRR_{ij}+\frac{1}{2}(\nabla_{i}\theta_{j}+\nabla_{j}\theta_{i})=
\big((n-2)(K-(h^\prime)^2)-hh^{\prime\prime}+hh^\prime\theta_t\big)g_{ij}^K\,.
\end{align*}
It is possible to see that, by construction, actually $\theta_i=0$ for
every $i$ at every point, but it is easier to consider directly
the 1--form $\sigma=\theta_t dt$ on $M^n$ and checking that it also
satisfies these three equations as $\theta$, hence the structural
equation
\begin{equation*}
\RRR_{\gamma\beta}+\frac{1}{2}(\nabla_\gamma\sigma_\beta+\nabla_\beta\sigma_\gamma)
=\frac{\alpha}{n} g_{\gamma\beta}\,.
\end{equation*}
It is now immediate to see that, $d\sigma_{it}=\nabla_{i}\sigma_{t}-\nabla_{t}\sigma_{i}=0$ and
$d\sigma_{ij}=\nabla_i\sigma_j-\nabla_{j}\sigma_{i}=0$, so the form
$\sigma$ is closed and being $M^n$ simply connected, there exists a
smooth function $f:M\to\RR$ such that $df=\sigma$, thus
\begin{equation*}
\RRR_{\gamma\beta}+\nabla^2_{\gamma\beta}f=\frac{\alpha}{n} g_{\gamma\beta}\,,
\end{equation*}
that is, the soliton is a gradient soliton.\\
It is also immediate to see that the function $f$ depends only on
$t\in I$.
\end{proof}

In the expanding, noncompact case (in the compact
case the soliton can be only a quotient of the hyperbolic space 
$\mathbb{H}^n$), if the Ricci tensor is nonnegative and 
$(M^n,g)$ is a gradient soliton, then either it is a warped product 
of $\SS^{n-1}$ (and  $M^n=\RR^n$) or it is the product of $\RR$ with a constant curvature
  space, but this last case is possible only if the soliton is the
  Gaussian expanding Ricci soliton, $\alpha<0$, on the flat $\RR^n$.

For a discussion of the expanding Ricci solitons which are
warped products of $\SS^{n-1}$ see~\cite[Chapter~1, Section~5]{chowbookI}, where
the authors compute, for instance, an example with positive Ricci
tensor (analogous to the Bryant soliton).\\
To our knowledge, the complete classification of 
complete, expanding, gradient, LCF Ricci solitons is an open problem,
even if they are rotationally symmetric.

\section{Singularities of Ricci Flow with Bounded
Weyl Tensor}

Let $(M^n,g(t))$ be a Ricci flow with $M^n$ compact on the maximal
interval $[0,T)$, with $T<+\infty$. Hamilton proved that
$$
\max_M|\Rm|(\cdot,t)\rightarrow\infty
$$
as $t\to T$.
 
We say that the solution has a Type I singularity if
$$
\max_{M\times [0,T)} (T-t) |\Rm|(p,t) <+\infty\,,
$$
otherwise we say that the solution develops a Type IIa singularity.\\

By Hamilton's procedure in~\cite{hamilton9}, one can choose a sequence of points 
$p_i\in M^n$ and times $t_i\uparrow T$ such
that, dilating the flow around these points in space and time, such 
sequence of rescaled Ricci flows (using Hamilton--Cheeger--Gromov
compactness theorem in~\cite{hamilton6} and Perelman's injectivity
radius estimate in~\cite{perel1}) converges to a complete maximal
Ricci flow $(M_\infty,g_\infty(t))$ in an interval $t\in(-\infty,b)$
where $0<b\leq+\infty$.\\
Moreover, in the case of a Type I singularity, we have $0<b<+\infty$, 
$|\Rm_\infty|(p_\infty,0)=1$ for some point $p_\infty\in M_\infty$ and
$|\Rm_\infty|(p,t)\leq 1$ for every $t\leq 0$ and $p\in M_\infty$.\\
In the case of a Type IIa singularity, $b=+\infty$,
$|\Rm_\infty|(p_\infty,0)=1$ for some point $p_\infty\in M_\infty$ and
$|\Rm_\infty|(p,t)\leq 1$ for every $t\in\RR$ and $p\in M_\infty$.

These ancient limit flows were called by Hamilton {\em singularity models}. We
want now to discuss them in the special case of a Ricci flow with uniformly bounded Weyl
tensor (or with a blow up rate of the Weyl tensor 
which is of lower order than the one of
the Ricci tensor). The Ricci flow under this condition is investigated also in~\cite{machen}.\\
Clearly, any limit flow consists of LCF manifolds, hence, by
Corollary~\ref{princ} and the cited results of Chen~\cite{chen2} and
Zhang~\cite{zhang} at every time and every point the manifold has nonnegative
curvature operator and either the Ricci tensor is proportional to the
metric or it has an eigenvalue of multiplicity $(n-1)$ and another of
multiplicity 1.

We follow now the argument in the proof of Theorem~1.1
in~\cite{caochen}.\\
We recall the following splitting result (see~\cite[Chapter~7,
Section~3]{chowluni}) which is a consequence of Hamilton's strong
maximum principle for systems in~\cite{hamilton2}.

\begin{teo}\label{split} Let $(M^n,g(t))$, $t\in(0,T)$ be a simply
  connected complete Ricci flow with nonnegative curvature
  operator. Then, for every $t\in(0,T)$ we have that $(M^n,g(t))$ is
  isometric to the product of the following factors,
\begin{enumerate}
  \item the Euclidean space,
  \item an irreducible nonflat compact Einstein symmetric space with
    nonnegative curvature operator and positive scalar curvature,
  \item a complete Riemannian manifold with positive curvature operator,
  \item a complete K\"ahler manifold with positive curvature operator on real $(1,1)$--forms. 
 \end{enumerate}
\end{teo}

Since we are in the LCF case, every Einstein factor above must be a
sphere (the scalar curvature is positive). The K\"ahler factors can be
excluded as the following relation holds for K\"ahler manifolds of
complex dimension $m>1$ at every point (see~\cite[Proposition~2.68]{besse})
$$
\vert\WWW\vert^2\geq \frac{3(m-1)}{m(m+1)(2m-1)}\RRR^2\,.
$$
Thus, any K\"ahler factor would have zero scalar curvature, hence would be
flat. Finally, by the structure of the Ricci tensor and the fact that
these limit flows are nonflat, it is easy to see that
only a single Euclidean factor of dimension one is admissible, moreover, in
this case there is only another factor $\SS^{n-1}$.\\
In conclusion, passing to the universal cover, the possible limit
flows are quotients of $\RR\times\SS^{n-1}$ or have a positive curvature
operator.

\begin{prop}[LCF Type I singularity models] Let $(M^{n},g(t))$, for 
$t\in[0,T)$, be a compact smooth solution to the Ricci flow with
  uniformly bounded Weyl tensor.\\
If $g(t)$ develops a Type I singularity, then there are two possibilities:
\begin{enumerate}
  \item $M^{n}$ is diffeomorphic to a quotient of $\SS^{n}$ and the
solution to the normalized Ricci flow converges to a constant positive curvature
metric.\\
In this case the singularity model must be a shrinking compact
Ricci soliton by a result of Sesum~\cite{sesum2}, hence by the analysis in
the previous section, a quotient of $\SS^{n}$ (this also follows by
the work of B\"ohm and Wilking~\cite{bohmwilk}).
\item There exists a sequence of rescalings which converges to the
  flow of a quotient of $\RR\times\SS^{n-1}$.
\end{enumerate}
\end{prop}
\begin{proof} By the previous discussion, either the curvature operator is
  positive at every time or the limit flow is a quotient of
  $\RR\times\SS^{n-1}$.\\
Hence, we assume that every manifold in the limit flow has 
positive curvature operator. The family of metrics 
$g_{\infty}(t)$ is a complete, nonflat, LCF, 
ancient solution with uniformly bounded positive curvature operator
which is $k$--non collapsed at all scales (hence a $k$--solution in
the sense of~\cite{perel1}). By a result of Perelman
in~\cite{perel1}, we can find a sequence of times $t_i\searrow-\infty$
such that a sequence of suitable dilations of $g_{\infty}(t_i)$
converges to a nonflat, gradient, shrinking, LCF Ricci soliton. Hence, we can
find an analogous sequence for the original flow. By the 
classification in the previous section, the thesis of the proposition
follows.
\end{proof}

\begin{rem}
Notice that in case (2) we are not claiming that every 
Type I singularity model is a gradient shrinking Ricci soliton.\\
This problem is open also in the LCF situation.
\end{rem}

\begin{prop}[LCF Type IIa singularity models] Let $(M^{n},g(t))$, for 
$t\in[0,T)$, be a compact smooth solution to the Ricci flow with
  uniformly bounded Weyl tensor. If the flow develops a Type IIa
  singularity, then there exists a sequence of dilations which
  converges to the Bryant soliton. 
\end{prop}
\begin{proof} As we said, if the curvature operator gets some zero eigenvalue, 
the limit flow is a quotient of $\RR\times\SS^{n-1}$ which cannot be a
steady soliton as it is not eternal. Hence, the curvature operator is
positive.\\
By Hamilton's work~\cite{hamilton13}, any Type IIa
  singularity model with nonnegative curvature operator and positive
  Ricci tensor is a steady, nonflat, gradient Ricci soliton. 
Since in our case such soliton is also 
LCF, by the analysis of the previous section, it must be the Bryant
  soliton.
\end{proof}

\bigskip

\begin{ackn} We thank Peter Petersen and G\'erard Besson for several 
valuable suggestions.\\
We also wish to thank 
Fabrizio Bracci, Alessandro Cameli and Paolo Dell'Anna for several
interesting comments on earlier versions of the paper.\\
The authors are partially supported by the Italian project FIRB--IDEAS
``Analysis and Beyond''.\\
The second author is partially supported by the Italian GNAMPA (INdAM) Group.
\end{ackn}

\bibliographystyle{amsplain}
\bibliography{LCFSolitons}

\end{document}